\documentclass[12pt]{article}
\usepackage{amssymb}
\usepackage{amsmath,amsthm}

\providecommand{\U}[1]{\protect\rule{.1in}{.1in}}
\makeatletter
\def\@seccntformat#1{\csname the#1\endcsname.\quad}
\makeatother

\begin{document}

\title{{\Large \textbf{Ultrametric Cantor sets and Origin of Anomalous Diffusion   }}}
\author{ Dhurjati Prasad Datta\thanks{
Corresponding author, email:dp${_-}$datta@yahoo.com}, Santanu Raut\thanks{%
{Chhatgurihati Seva Bhavan Sikshayatan High Scool, New Town,
Coochbehar-736101, India, email}: raut\_santau@yahoo.com} and Anuja Roy
Choudhuri\thanks{
Ananda Chandra College, Jalpaiguri-735101, India, email: anujaraychaudhuri@ymail.com} \\
Department of Mathematics, University of North Bengal \\
Siliguri,West Bengal, Pin: 734013, India }
\date{}
\maketitle

\begin{abstract}
The anomalous mean square fluctuations are shown to arise naturally from the  ordinary diffusion equation interpreted scale invariantly in a formalism endowing real numbers with a nonarchimedean multiplicative structure. A variable $t$ approaching 0 linearly in the ordinary analysis is shown to enjoy instead a sublinear $t\log t^{-1}$ flow in the presence of this scale invariant structure. Diffusion on an ultrametric Cantor set is also generically  subdiffusive with the above sublinear mean square deviation. The present study seems to offer a new interpretation of a possible emergence of complex patterns from an apparently simple system.
\end{abstract}

\begin{center}


\end{center}

\baselineskip=16.5pt

\newpage

\section{Introduction}
 Anomalous diffusion is known to occur ubiquitously in diverse 
 complex systems enjoying fine ``structure with variations" \cite{compx} such as in disordered or fractal media, thus developing fat tailed, broad distributions and/or long range spatio-temporal correlations \cite{sub, sup}. The hallmark of such a diffusion process is the occurrence of an anomalous law for the mean square displacement (fluctuation), viz., $ \langle \Delta x^2(t)\rangle=t^{\nu}$, with $\nu\neq 1$. Sub-diffusive ($\nu < 1$) behaviour is usually predominant  in disordered systems,  for instance, in spin glasses, amorphous semiconductors, lipid bilayers, living cells, transport in fractal sets and many others where a broad (fat tailed) distribution for local trapping times of the diffusive test particle gradually build up \cite{sub}. Super-diffusion ($\nu > 1$), on the other hand, may arise from long range correlations in velocity fields of turbulent flows, Levy flights and so on \cite{sup}. The emergence of nonlinear growth of the mean square fluctuation clearly reflects a possible violation of the gaussian central limit theorem in the underlying random walk processes \cite{comment}. An important, still unsolved, problem is to look for a generic (universal) mechanism for the emergence of anomalous diffusion in such diverse phenomena. In the present paper, we offer one such  dynamical principle that might be  at play at the heart of complex systems. We show that the anomalous mean square fluctuations can arise naturally from the  ordinary diffusion equation interpreted scale invariantly in a formalism endowing real numbers with a nonarchimedean multiplicative structure. A variable $t$ approaching 0 linearly in the ordinary analysis is shown to enjoy  a sublinear $t\log t^{-1}$ flow in the presence of this scale invariant structure. Diffusion on an ultrametric Cantor set is also generically  subdiffusive with the above seemingly universal sublinear mean square deviation. The present study seems to offer a new interpretation of a possible emergence of complex patterns from an apparently simple system. This, in turn, appears to realize the philosophical musings, `` Nature can produce complex structures even in simple situations, and can obey simple laws even  in complex situation \cite{compx}."

There is already a vast body of studies on anomalous diffusion and its origin that are available in literature \cite{sub,sup}.  Even at this back ground,  the present investigation aims at offering a potentially  new insight into the actual mechanism of the dynamics of an anomalous motion. As
it is well known, there are actually two distinct types of motion observed in
Nature: smooth, regular motion, like the Newtonian (two body) planetary motion, and
random, highly irregular motion, as in the Brownian motion of a fine pollen
particle in a liquid at rest \cite{sup}. A smooth motion, at least on a moderate time scale, is expected to be predictable and so are deterministic in nature, when Brownian type motion requires statistical ( stochastic) methods. Of course, the deterministic chaos falls in between, and several authors, for instance, Ref. \cite{dyn}, discussed the problem of offering  a dynamical interpretation  for the Brownian like motion based on the deterministic Hamiltonian models in the phase space.  At a more elementary level, on the other hand, the classical Brownian motion can, in fact, be
considered to enjoy a bit of a deterministic flavour as the relevant
gaussian probability distribution (transition probability/ propagator) is
known to follow the linear homogeneous diffusion equation 
\begin{equation}  \label{diff1}
{\frac{\partial  W}{\partial t}}= {\frac{\partial^2  W}{\partial x^2}}, \ -\infty<x<\infty, \ t>0
\end{equation}
\noindent where we choose the diffusion constant to be unity. Accordingly,
variations of macroscopic variables such as the number of diffusive
particles, concentration and similar other quantities are governed
effectively by a smooth deterministic law. Mathematically, Brownian motion
is a process realized in a homogeneous smooth manifold where Taylor's
expansion and other relevant analytic function theoretic resources are
available. Moreover, the probability density is also a smooth function
having finite mean and variance, so that the universality of the central
limit theorem drives the force law to be smooth (for another explanation, see \cite{dyn}).

On the otherhand, if the underlying space of diffusion has a manifestly
disordered, fractal structure, the above gaussian, and/or the function
theoretic smoothness is generally lost (see for instance \cite{barlow, fgz, pb, bukh, karo, sub}). All these studies tried to offer  precise mathematical justifications leading
to anomalous mean square variations in such fractal spaces \cite{comment1}. It follows also
that the emergence of a {\em smooth} effective deterministic law (in the sense of  differentiable functions in a Euclidean space) at the macroscopic
scale is generally lost because of the inherent loss of smoothness in the random
process \cite{comment2}.

Over the past few years we have been developing a formalism of a scale
invariant analysis \cite{dp1,dp2} aiming at integrating the framework of the standard
analysis on, for instance, the Euclidean spaces and those on a fractal like
space. Another motivation was to investigate {\em if a seemingly smooth
deterministic evolution, even in the absence of any external influences,
may lead naturally to a complex pattern over an asymptotically long (or short) time
scale.} In the following subsection, we present a short summary of the basic ingradients of our approach \cite{dp1,dp2,dp3,dp4} that will be necessary for the subsequent analysis. However, most of the derivations are presented first time here offering better insights as compared to the original derivations. In the latter two sections, we present new results on Cantor sets and diffusion on a scale invariant real line and Cantor sets respectively.
It will become clear that the presence of {\em dynamically} generated ultrametric Cantor sets at infinitesimally small finer scales in the real line  would force a linear diffusive process to evolve anomalously over infinitely long time scales, that are available naturally in the scale invariant formalism.

\subsection{ Formalism}

We construct a scale invarinat, nonarchimedean extension of the real number system $R$ accommodating the so called generalized solutions of the scale free Cauchy problem
\begin{equation}\label{ssim}
t{\frac{dx}{dt}}=x, \ x(1)=1
\end{equation}
\noindent in the form 
\begin{equation}\label{ns}
X_{\pm}(t)=t\tilde X_{\pm}(\tilde t)  
\end{equation}

\noindent where $\tilde X_{\pm}=\tilde t^{\mp v(\tilde t)}, \ \tilde t=t/\epsilon$, and $v$ satisfies the self similar replica of an inverted (\ref{ssim}), viz,

\begin{equation}\label{ssim2}
\tau{\frac{dv}{d\tau}}=-v
\end{equation}
\noindent in the logarithmic variable $\tau=\log \tilde t$, for an arbitrarily small $\epsilon$. Notice that $X_{\pm}(1)=1$ as one must set $\epsilon =1$ and $v$ is finite. Clearly, $v\equiv 0$ in the ordinary analysis on the real line. However, in an extended nonarchimedean model $\bf R$ of the real number system $ R$ accommodating {\em nontrivial infinitesimals} in the neighbourhood of zero, this class of solutions becomes meaningful (for a justification see \cite{dp2}). Notice that the ansatz (\ref{ns}) tells that each element $t\in R$ splits into pairs of infinitesimally close neighbours belonging respectively to the right and left neighbourhood of $t\in \bf R$.  The scale free component $\tilde X$ (for simplicity of notations, we drop temporarily subscripts $\pm$) of the generalized solution  now is realized as a locally constant function (LCF) defined by 
\begin{equation}\label{lcf}
t{\frac{d\tilde X}{dt}}=0
\end{equation}
\noindent As a consequence, $\tilde t=\tilde t(t)\in (0,1)$ may, in fact, be {\em any} continuously first differentiable function of $t$, which is a consequence of the {\em reparametrisation invariance} of a locally constant function in an ultrametric space (this aspect of a LCF is discussed in more detail in \cite{dp3}).

Now to construct a self consistent ultrametric extension suitable for our purpose, let us introduce scale invariant infinitesimals via a more refined evaluation of the limit $t\rightarrow 0^+$ in ${R}$. Notice that as $t\rightarrow 0^+$, there exists $\epsilon>0$ such that $0<\epsilon<t$ and one may like to identify zero (0) with the closed interval $I_{\epsilon}=[-\epsilon,\epsilon]$ at the {\em scale} (i.e. accuracy level in a computation) $\epsilon$. Ordinarily, $I_{\epsilon}$ is a connected line segment, which shrinks to the singleton $\{0\}$ as $\epsilon\rightarrow 0^+$, so as to reproduce the infinitely accurate, {\em exactly determinable} ordinary real numbers. 

Let us now present yet another {\em nontrivial mode} realizing the limiting motion $t\rightarrow 0^+$. This mode is nonlinear, as it is defined via {\em inversions}, rather than simply by the linear translations, that is available {\em uniquely} in the standard analysis  on $R$. In the presence of {\em infinitesimals}, the present nonlinear mode may, nevertheless, assume significance. In the following we give a construction of an {\em explicitly defined infinitesimal}, without requiring the model theoretic set up of Robinson's nonstandard analysis. {\em Our formalism may be considered to be a new, independent realization of inifinitesimals residing originally in an ultrametric space, but, nevertheless, inducing {\large real valued ``infinitesimal corrections"} to the ordinary real variable} \cite{comment3}. A real number, as it were, assumes a {\em deformed value} because of a possible nontrivial {\em motion} close to 0. Fix $\epsilon=\epsilon_0$ and let $C_{\epsilon_0}\subset [0, \epsilon_0]\equiv I_{\epsilon_0}^+$ be a Cantor set defined by an IFS of the form 

\begin{equation}\label{ifs}
f_1(t)=\lambda t, \ f_2(t)=\lambda t-(\lambda/\epsilon_0 -1)\epsilon_0
\end{equation}
\noindent where $\lambda=\beta\epsilon_0, \ 0<\beta<1$ and $\alpha+2\beta=1$. Thus, at the first iteration an open interval $O_{11}$ of size $\alpha\epsilon_0$ is removed from the interval $I_{\epsilon_0}^+$, at the second iteration 2 open intervals $O_{21}$ and $O_{22}$  each of  size $\alpha\epsilon_0(\beta)$ are removed and so on, so that a family of gaps $O_{ij}$ of sizes $\alpha\epsilon_0(\beta)^{i-1}, \ j=1,\ldots, 2^{i-1}$  are removed in subsequent iterations from each of the closed subintervals $I_{ij}, \ j=1,2,\ldots, 2^{i} $ of $I_{\epsilon_0}^+$. Consequently, $C_{\epsilon_0}=I_{\epsilon_0}^+ \ -\ \underset{i}{\cup}O_{ij}=\underset{i}{\cap}\underset{j}{\cup}I_{ij}$. Notice that the total length removed is $\sum \alpha\epsilon_0(2\beta)^{i-1}=\epsilon_0$, so that $m(C_{\epsilon_0})=0$.

Next, consider $\tilde I_{N}=[0, \beta^N]$ and let $N=n+r$ and $N\rightarrow \infty$ as $n\rightarrow \infty$ for a fixed $r\geq 0$. Choose the scale $\epsilon=\alpha \beta^n\epsilon_0$ and define $\tilde t_r\in [0, \alpha\beta^N\epsilon_0]$ a ( positive ) {\em relative infinitesimal} (relative to the scale $\epsilon$) provided it also satisfies the {\em inversion} rule
\begin{equation}\label{inv}
\tilde t_r/\epsilon\propto \epsilon/t, \ 0<\tilde t_r<\epsilon<t
\end{equation}
\noindent In the limit $\epsilon\rightarrow 0$, relative infinitesimals $\tilde t_r$, of course, vanish identically. However, the corresponding scale invariant infinitesimals ${\tilde \eta}_r=\tilde t_r/\epsilon, \ \epsilon \rightarrow 0$ may nevertheless be nontrivial. As a consequence, the relative infinitesimals may be awarded a new {\em scale invariant absolute value} \cite{dp1,dp2,dp3}
\begin{equation}\label{value}
v(\tilde t_r)\equiv |\tilde t_r|_u=\underset{\epsilon\rightarrow 0}{\lim}\log_{\epsilon^{-1}}\epsilon/|\tilde t_r|, \  \tilde t_r\neq 0
\end{equation}
\noindent and $|0|_u =0$. The set of infinitesimals are uncountable, and the above norm satisfies the stronger triangle inequality $|a+b|_u\leq {\rm max}\{|a|_u,|b|_u\}$. Accordingly, the zero set ${\bf 0}= \{0, \pm\epsilon {\tilde \eta}_r| \ {\tilde \eta}_r\in (0,\beta^r), \ r=0,1,2,\ldots, \ \epsilon \rightarrow 0^+\}$ may be said to acquire {\em dynamically} the structure of a Cantor like ultrametric space, for each $\beta\in (0,1/2)$ (so as to satisfy the open set condition \cite {dp3}). The set $\bf 0$ indeed is realized as a nested circles $S_r: \{\tilde t| \ |\tilde t_r|_u=\alpha_r\}$, in the ultrametric norm, when we order, without any loss of generality, $\alpha_0>\alpha_1>\ldots$. For definiteness, we call 0 the hard (or stiff) zero, whereas elements of the asymptotic form $\pm\epsilon {\tilde \eta}_r, \ \epsilon \rightarrow 0^+$ are designated as soft (or dynamic) zeros.

To give a more prictorial representation of infinitesimals, let us recall that in a computational problem, a number $x=1$, for example, is represented upto a finite acuracy; i.e., upto to a scale  $\epsilon_0$, say. Then the numbers in the interval $(1-\epsilon_0,1+\epsilon_0)$ are {\em computationally unobservable} and identified, as a whole, as the number 1. The above nontrivial construction now tells that the points  in that computationally inaccessible  limiting interval are {\em aligned dynamically} as nonintersecting clopen balls of a Cantor set $C$ endowed with the nontrivial ultrametric value (\ref{value}), thereby extending the ordinary set $R$ to an infinite dimensional, scale free, nonarchimedean space $\bf R$ accommodating {\em dynamically active} infinitesimals and infinities [c.f. Sec. 4]. A nonzero real variable $t (>0$ say,) in $R$ and approaching 0 now gets a deformed structure in $\bf R$ of the form $T(t)=t.t^{-v(\tilde t(t))}$, thus reflecting a nontrivial effect of dynamic infinitesimals over the structure of the real number system $R$.

 Next, we see that the above definition (\ref{value}) is well defined and is not empty. For, let us assign $\phi(\tilde t)=v(\tilde t_r)\log {\tilde t_r}^{-1}$ a constant valuefor all ${\tilde\eta}_r\in (0, \beta^r)$. This constitues a useful class of ultrametric norms when $\phi(\tilde t)$ is identified with a {\em Cantor function} \cite{dp3}. Indeed, for each arbitrarily small, but nonzero, $\epsilon\neq 0$, relative infinitesimals are elements of gaps of a Cantor set of the type $C_{\epsilon}$, so that $\phi(\tilde t_r)$ is assigned constant values, in a piecewise sense,  every where in $[0,\epsilon]$ so that $d\phi/dt=0$ a.e in $[0,\epsilon]$. Indeed, the constant values can always be assigned continuously  so as to make it continuous and non-decreasing on $[0,\epsilon]$.  However, at the end points of the gaps (that is, at the elements of $C_{\epsilon}$), $\phi$ need not be differentiable in the usual sense. Hence, $\phi$ is uniquely realized as the Cantor function associated with $C_{\epsilon}$.

In the present approach, however, $\phi$ is not only continuous, but also satisfies $\frac{d\phi}{dt}=0$ every where in $I_{\epsilon}^+$ (for details, see, for instance, \cite{dp2,dp3}). The usual derivative discontinuity is removed, since in the vicinity of $t_0\in C_{\epsilon}$, transition between the end  points of the gap containing $t_0$ is now accomplished smoothly by inversions of the form (\ref{inv}).

To summarise, the nontrivial infinitesimals are represented scale invariantly by the asymptotic formula of the form  $\tilde X_{\pm}=X_{\pm}/x=x^{\mp v(\tilde x)}$ as $x\rightarrow 0$ (from now on, $t$ and $x$, as usual, denote time and space variables). Infinitesimals may therefore be interpreted as those numbers which approach 0 at a slower (nonlinear ) rate: $\log \tilde X_-^{-1}=\log \tilde X_+=v(\tilde x)\log \tilde x^{-1}$ as $x\rightarrow 0$. The solution (\ref{ns}) then provides a nontrivial representation for the {fattened} real
number $x$, denoted as $\mathbf{x}$. The infinitesimals in a more
conventional nonstandard model \cite{robin} of $R$ are, however, inactive or passive in
the present sense. The dynamical infinitesimals are already shown to have
nontrivial influences in the asymptotic estimates of number theory and other
areas of real analysis \cite{dp2,dp3}. Further studies will be reported elsewhere. For our
latter reference, let us note here that one does enjoy a large body of resources of
dynamical infinitesimals residing in a nontrivial neighbourhood of zero (0).
Such an infinitesimal may be assumed to live in an infinitesimal
copy of an arbitrarily assigned  Cantor set. One has to make a \emph{choice} of an appropriate
Cantor set depending on the problem in hand. In the following we analyse the
implications of the above nonarchimedean structure on a diffusion process on
a measure zero Cantor set and also on $\bf R$.

\section{Cantor set: New Results}

The above ultrametric structure can easily be realized on a Cantor
set (a compact and perfect subset of $R$). Recently, it is shown that 
the ultrametric so defined is both
metrically and topologically inequivalent to the usual ultrametric that a
Cantor set  carries naturally \cite{dp3}. For a
point $x_0$ in a Cantor set $C\subset I=[a, b]$, the representation (\ref{ns}%
) now gives rise to a scale invariant ultrametric extension $\tilde
X_{\pm}=X_{\pm}/x_0=x_0^{\mp v(\tilde x)}$ where the transition between two infinitesimally
close scale invariant neighbours is supposed to be mediated by inversions of the form $
\tilde X\rightarrow {\tilde X}^{-h}$ for a real $h$, which determines the
jump size. Notice that $\tilde X$ (and equivalently, $v(\tilde x)$) is a locally
constant Cantor like function and solves ${\frac{d\tilde X}{dx}}=0$
everywhere in $I$. The ordinary discontinuity of a Cantor function at an $
x_0\in C$ is removed, since in the present ultrametric extension, the point $
x_0$ in $C$ is replaced by an \emph{inverted Cantor set} which is the
closure of \emph{gaps} of an infinitesimal Cantor set $C_i$ that is assumed
to be the residence Cantor set for the relevant infinitesimals $\tilde x$ living in the extended neighbourhood $\bf 0$ of 0. The gaps of $C_i$ constitute a disjoint family of connected clopen intervals (represented
in a scale invariant manner) over each of which scale invariant equations of
the form (\ref{ssim}) are well defined \cite{dp3}. Consequently, the valuation $v(\tilde x)$, redefined slightly in the modified form 

\begin{equation}\label{value1}
(x/x_0)^{\tilde {v}(x)}=x_0^{v(\tilde x)} 
\end{equation}

\noindent (that is, $\tilde v(x)/v( \tilde x)=\log x_0/\log (x/x_0)$, 
exposing the relative variation of $\tilde v$ over $v$), $x$ assuming values 
from the gaps in the neighbourhood
of $x_0$, is realized as a smooth function defined recursively in a scale
invariant way by the equation

\begin{equation}  \label{valueq}
{\frac{d\tilde v(x)}{d \xi}}=-\tilde v(x)
\end{equation}

\noindent where $\xi=\log\log (x/x_0), \ x\in C_i$. Recall that $\tilde x$ resides in the gaps of nontrivial neighbourhood of 0 instead. As a consequence, $\tilde v$ may 
be written as $\tilde v(x)=(\log (x/x_0)^k)^{-1}$, where $k$ may be allowed to assume values from a set of scale factors related to that of the Cantor set. This form is clearly consistent with (\ref{value1}). Assuming $x$ is drawn from a specific gap of a given size, the same, written more effectively as $\tilde v(x)/v(\tilde x)=(\log_{x_0}(x/x_0))^{-1}$, yields, in the limit of  vanishingly small gaps (i.e., as $x\rightarrow x_0$ and vice versa), the limiting value $\tilde v_0(x)/v_0( \tilde x)=1/s$, since {\em $ \lim \log_{x_0}(x/x_0)=s$ equals the finite Hausdorff dimension of the Cantor set $C_i$.} Let us first note that if one replaces the Cantor set by a segment of a line of the form $(0, \delta)$, then $x/x_0=1$, in that limit ($\delta\rightarrow 0$) gives $s=0$, which is consistent with the fact that the line segment reduces to a point, viz. 0. In the general case, $x/x_0\propto N$, the number of clopen balls that covers the fattened gap of the form $(x,x_0)\subset (0, \delta)$ (size of balls are determined by the gap). Letting $x_0\rightarrow 0$ (as the relevant scale factors $\beta^n$), the above limit therefore mimics the box dimension, which also equals the Hausdorff dimension of the  Cantor set concerned.
The topological inequivalence of the present ultrametric arises from the possible dichotomy in the choice of $C_i$. 

Let us remark that the ordinary limit $\epsilon\rightarrow 0$ in $R$ is altered because of scale invariant dynamic infinitesimals $\log \tilde x/\epsilon \approx v(\tilde x)\log \epsilon^{-1}\approx \epsilon\log \epsilon^{-1}$, when $\tilde x$ is considered to lie on a fattened (connected ) gap, so that the ultrametric valuation  may be assumed to coincide with the usual (Euclidean) value viz. $v(\tilde x)\approx \epsilon$. However, assuming $\tilde\epsilon (=\beta^n, \ n \rightarrow \infty)$ to be an infinitesimal scale of the Cantor set concerned, we also have $\log \tilde x/\tilde\epsilon \approx \tilde\epsilon^s\log \tilde\epsilon^{-1}$, since the valuation is identified with the associated Cantor function $\phi (\tilde x)\approx \tilde \epsilon^s\approx \epsilon$, $s$ being, as usual, the corresponding Hausdorff dimension. Reverting back to the ordinary scale $\epsilon$ (and keeping in view the associated scale invariance), this scaling can be identified with $\epsilon^{\tilde s}\approx O(1)\epsilon\log \epsilon^{-1}$, for an $\tilde s$ given by $\tilde s\approx 1-\frac{\log\log \epsilon^{-1}}{\log \epsilon^{-1}}$. {\em As a consequence, in the presence of an ultrametric space, the ordinary limit $\epsilon \rightarrow 0$ is replaced by the sublinear limit 

\begin{equation}\label{sublin}
\epsilon^{\tilde s}=\epsilon\log \epsilon^{-1}\rightarrow 0, 
\end{equation}

\noindent $ 0<\tilde s<1$, as $\epsilon\rightarrow 0$}. This is one of the main results of this paper. A real variable $t$ in $R$ approaching to (or flowing out from ) 0 will experience this scale invariant sublinear behaviour in an incredibly small neighbourhood of $\bf 0$ in $\bf R$ and should have a deep significance in number theory and other areas \cite{dp2}. Similar behaviour is also reported recently  in the context of diffusion in an ultrametric Cantor set in a noncommutative space \cite{pb}. 
Finally, the ultrametric induced by the valuation $v(\tilde x)$
coincides with the natural ultrametric only when the scaling properties of $
C_i$ coincide with that of $C$. In this paper we adhere to the latter possibility.

Suppose the original Cantor set $C$ and the infinitesimal Cantor set $C_i$
have Hausdorff dimensions $s$ and $s^{\prime}$ respectively. Any
point $\mathbf{x}$ of the fattened set $\mathbf{C}= C \ + \ C_i$ is given as $%
\mathbf{x}=x \ + \ \tilde x, \ x\in C, \ \tilde x\in C_i$. It is well known
that $\mathbf{C}=I$, for a.e. $s^{\prime}$, for a given $s$ \cite{solo}.
Accordingly, it follows that  given a Lesbegue measure zero Cantor set $C
$, the above smooth differentiable structure is a.s realized on $\mathbf{C}$%
, which is nothing but $I$, though in an appropriate (scale free)
logarithmic variable.

To understand more clearly the above smooth scale invariant structure let us
consider the classical middle third Cantor set $C_{1/3}$ with scale factors $%
\tilde \epsilon=3^{-n}$. A point $x_0=3^{-n}\sum a_i3^{-i}, \ a_i\in \{0,2\}$ of $C_{1/3}
$ is raised to the scale free $\mathbf{x}$ which is a variable living in a
family of fattened gaps, attached and structured hierarchically at the point $x_0$ (or
equivalently, by scale invariance, at 1), over each of which scale free
equations of the form (\ref{valueq}) are valid. The infinitesimals are elements of the gaps ``closest" to 0, viz. the open intervals $3^{-n-m}(1,2)$, in the limit $n\rightarrow \infty, $ for a fixed  $0<m<n$, which are assigned  nontrivial values akin to the Cantor function $v(\tilde x)=i3^{-ms}, \ i=1,2,\ldots 2^m-1$ \cite{dp1}.  
Over each of the finite size  gaps, on the otherhand,  the
valuation $v(x)$ is awarded as $v(x)=3^{-sn}$. Both these valuations are not only continuous but also smooth  since the corresponding Cantor function is realized as a
smooth function via the logarithmic ansatz for a substitution of the
form $3^{n}\Delta x_n=3^{n}(x-x_n)\rightarrow n\log {\frac{\mathbf{x}}{x_n}}$, as $
n\rightarrow \infty$, thereby removing the derivative discontinuity at the
points of scale changes \cite{dp1}, so that ${\frac{dv(x)}{dx}}=0$, every where on the
Cantor set concerned. Notice that gaps scale as $\epsilon=2^{-n}$ (recall the binary representation for points on a connected segment of the real line) when a closed interval containing  points like $x_0\in C_{1/3}$  scales as $3^{-n}$, so that the Hausdorff dimension is $s=\log_3 2$.
By (\ref{value1}), the variability of the valuation $\tilde v(x)$ in the limit of vanishing gap sizes is obtained as  $\tilde v_0(x)\propto 3^{-sn}s^{-1}, n\rightarrow \infty$.

 Next, before determining the incremental \emph{measure}, denoted $d_j\tilde X$,
of smooth self similar jump processes of (gap) ``size" (in the sense of a weight) $\epsilon \ (2^{-n},$ for $C_{1/3}$, say) in the neighbourhood of the
scale invariant 1, let us first recall that pure translations follow a linear law: $y=Tx=x+h$. The {\em instantaneous pure jumps} (of unit length close to the scale invariant 1), on the other hand, follow a {\em hyperbolic law}: $\tilde X\rightarrow Y=\tilde X^{-1} \ \Rightarrow \log Y+\log X=0$, which tells, in turn, that the corresponding translational increment, even in the log scale, is indeed zero. This actually is the case for the valuation (\ref{ns}) defined in terms of the locally constant Cantor function. The (manifestly scale invariant) multiplicative valuation $\tilde v(x) = \log_{x^{-1}}(X/x)$ in (\ref{ns}), however, gives the correct linear measure for a single jump relative to the point $x$ (for the above hyperbolic type jump, $v(x)=1$ relative to $x$ itself as the scale). The corresponding {\em multiplicative increment} is denoted as $\delta_j \tilde X=(x/x_0)^{\tilde v(x)}$. More importantly, this valuation is realized as a smooth measure and may be considered to contribute an independent component in the ordinary measure of $R$. Further, the total self similar jump mediated increments over a spectrum of gaps of various sizes of the forms $\epsilon_n=2^{-n}\epsilon_0$ in the neighbourhood of a (middle third) Cantor point $x_0$ (say) is now obtained as 

\begin{equation}\label{jump}
\Delta_j \tilde X= (x/x_0)^{s^{-1}2^{-m}\underset{n}\sum 2^{-n}}
\end{equation}
\noindent which in the limit $x\rightarrow x_0$, that is, $m\rightarrow \infty$ yields the {\em jump differential} $d_j\tilde X= (d\tilde x)^{s^{-1}}$, where $\tilde x = \lim (x/x_0)^{2^{-m}},$ is a deformed variable close to 1. Such a variable ($\neq 1$ exactly) exists because of a nontrivial g.l.b of gap sizes (another manifestation of the sublinear asymptotic). Clearly, the jump differential appears to relate to a fractional differential of order $s^{-1}$. Invoking further the inherent smoothness of the formalism we, however, treat this  fractional differential as the ordinary differential $d_j\tilde X=d \tilde x^{s^{-1}}$. One justification, at least huerisitically, is the following. The jump process, although realized here as smooth, carries an inherent random element \cite{comment2}. The jump measure must therefore be interpreted as an equality of moments $<d\tilde x>^{s^{-1}}\approx  d<\tilde x^{s^{-1}}>$. A possible distribution with analogous  moment constraints may be Poisson like but not exactly Poissonian, an example of which is considered in \cite{dp4}.
Incidentally, we note that the essential singularity in $s=0$ tells that
in the absence of inversion mediated jumps, the whole structure of gaps
collapses to a point (singleton set, devoid of any nontrivial infinitesimals). The divergence in the jump measure
then reflects  the ordinary nondifferentiable structure of the Cantor
set. On the otherhand, on any connected segment of $R$, $s=1$, and the jump measure reduces to the ordinary linear measure $dx$. To summarise, {\em  the significantly new insight that emerges from the above analysis is that an infinitesimal scale invariant increment on an ultrametric space must have the form $\tilde X=1+{\epsilon}^{1/s},\ \epsilon\rightarrow 0^+$ on a connected segment close to 0.} Recalling $\tilde\epsilon={\epsilon}^{1/s}$, the above infinitesimal jump increment $\tilde X= 1\pm\tilde\epsilon$ reduces to the usual increment on a Cantor set in the usual metric, but at the cost of the smooth structure.

\section{Diffusion}

Coming back to the diffusion equation (\ref{diff1}), let us next recall the
self similarity of the same. Writing $W(x,t)=t^{-1/2}w(z), \ z=\frac{x}{%
\sqrt t}$, the scaling function $w$ satisfies the first order ordinary
differential equation 
\begin{equation}  \label{ode}
{\frac{dw}{du}}=-w, \ u=z^2
\end{equation}

\noindent giving rise to the gaussian propagator $W(x,t)=At^{-1/2}e^{-{\frac{
x^2}{4 t}}}$ by a direct integration. 

Now, in the present analysis, the real variables $x$ and $t>0$ must be assumed to live in the corresponding scale invariant extensions $\bf R$ and $\bf R_+$ (set of nonnegative numbers), so that the scaling variable $z$ gets extended to a deformed variable $\tilde Z=z/z_0={\frac{\tilde X}{\sqrt{\tilde T}}}\in \bf 1$ in the extended neighbourhood of a point $(x_0, t_0)\in R\times R_+$. The scale invariant variables $\tilde X$ and $\tilde T$, in turn, belong to two Cantor sets $C_s$ and $C_t$ respectively with scale invariant measures $d_j \tilde X=d\tilde x^{\alpha}$ and $d_j\tilde T=d\tilde t^{\beta}$, 
where $\alpha$ and $\beta$ are the
respective inverse Hausdorff dimensions and $\tilde x$ and $\tilde t$ are two scale invariant variables near 1 of $R$. As a consequence, Eq(\ref{ode}), defined originally on $R$, now automatically  gets extended to one in  the new rescaling symmetric variable $\tilde Z$ living on a ultrametric Cantor set in the extended neighbourhood of every point $(x, t)\in R\times R_+$
and hence supports nontrivial solutions analogous to (\ref{ns}). Thus, integrating (\ref{ode}) in that class of new solutions, we get 
\begin{equation}  \label{nongauss}
W_C(\tilde X, \tilde T)=A {\tilde t}^{-\beta/2}e^{-\lambda({\frac{\tilde x^{2\alpha}}{ \tilde t^{\beta}}})^{1+\nu}}
\end{equation}
\noindent as a {\em stretched exponential} fat tailed propagator for a diffusive process (random walker) on a Cantor set and $\lambda$ is a constant depending on $z_0$. Clearly, $W_C$ satisfies the scale invariant diffusion equation 
\begin{equation}  \label{diff2}
{\frac{\partial  W}{\partial \tilde T}}= {\frac{\partial^2  W}{\partial \tilde X^2}}
\end{equation}
\noindent which is defined close to every point $x\in C_s$ and at any instant $t\in C_t$. Although the derivatives are evaluated with jump differentials, these are equivalent to the usual partial derivatives, but in the deformed (scaling) variables $x^{\alpha}$ and $t^{\beta}$ respectively. Further, the exponent $\nu$ in (\ref{nongauss}) is a valuation so that ${\tilde Z}^{\nu(\tilde Z)}$ is a locally
constant Cantor function that arises in connection with the residence Cantor
set for the variable $z$ (c.f. eq(\ref{ns})) \cite{note1}). Because of the  sublinear asymptotic increments of the form (\ref{sublin}), this equation is also considered to be valid on a connected line segment close to $t=0$, so that (\ref{diff2}) is also valid for $\tilde t\approx 0$. Next, we note that the
scale invariant factors of a real variable (viz. eq(\ref{ns})) become significant only for an asymptotically large time. Consequently, the above scale invariant solution (\ref{nongauss}) of the diffusion equation (\ref{diff1}) is expected to arise naturally in any diffusive process that persists over many longer time scales living in a set of the form $\bf R$. 
The  anomalous mean square deviation is  given generically as $\langle \Delta x^2(t)\rangle= { t}^{\beta/ \alpha}$, where $t\rightarrow 0^+$, for every scale invariant $x$ near 1. 

We conclude that a simple diffusion process if allowed to evolve over many longer and longer time scales as those available for natural processes will ultimately give away naturally to  a {\em stetched exponential } nongaussian distribution of the form (\ref{nongauss}) leading to an anomalous mean square fluctuation; reflecting, in turn, the universally present scale invariant numerical fluctuations \cite{volv}.

Diffusion on a Cantor set when examined in the framework of the
ordinary Newtonian time  (i.e., when $C_t$ reduces to the singleton $\{0\}$) is generally subdiffusive with exponent $\alpha^{-1}=s:0<s<1$, $s$ being the Hausdorff dimension of the diffusive medium (since $\beta =1$). The mean square deviation has the generic form  $\langle \Delta x^2(t)\rangle= {t}^s\approx t\log t^{-1}, \ t\rightarrow 0^+$, because of the sublinear asymptotic flow on a Cantor set. As already remarked, because of scale invariance this may also be interpreted as the asymptotic late time anomalous mean square fluctuation from the gaussian law even for the ordinary diffusion equation (\ref{diff1}).

For a fractal time process the subdiffusion occurs for $s<\tilde s$ and superdiffusion for $s>\tilde s$, $\tilde s$ being the Hausdorff dimension for underlying Cantor like set for the fractal time. However, for $s=\tilde s$, the gaussian like linear mean square variation may be observed even for a fractal time process \cite{comment}. Analogous smoothening in the asymptotic scaling of the eigen value counting function was also noticed by Freiberg \cite{fgz}.

\section{Final Remarks}
Before closing, let us note that the usual ``cut and delete" process realizing a Cantor set seems to give a misleading idea that the Lebesgue measure of a set, for instance, [0,1] becomes zero under a recursive application of a contraction mapping of the form $f(x)=x/3$. However, the full Lebesgue measure could actually be preserved  multiplicatively at every level: $1=3^{-n} \times q^{ns}$ for a given $q>0$, leading to a Cantor set of dimension $s= \frac{\log 3}{\log q}$, that remains attached to the limit point of the said contraction.  It is reasonable to imagine that the (static) points of the interval [0,1] dynamically adjust among themselves in a {\em scrambled } manner to accommodate around the limit point a Cantor set with fractal dimension $s$. The intrinsic motion of these dynamic numbers (in the logarithmic scale) could be visualized as an anomalous Brownian motion given by (\ref{diff2}). Nature seems to make use in plenty this seemingly universal principle of {\em making space out of nothing!}


\begin{thebibliography}{99}
\bibitem{compx} N. Goldenfeld and L. P. Kadanoff, Science, {\bf 28}, 87-89, (1999).
\bibitem{sub} R. Metzater and J Klafter, Phys. Reports. {\bf 339}, 1-77, (2000).
\bibitem{sup} A. A. Dubkov et al, Int. J. of Bifurcation and Chaos, {\bf 18}, 2649 - 2672 (2008).
\bibitem{comment} Complex systems with linear mean square fluctuation in Browninan like motion with nongaussian distribution are reported in literature recently,  B. Wang et al, PNAS, {\bf 206}, 15140-15146, (2009). See also J. K. E. Tunley, Phys. Rev. Lett, {\bf 33}, 1037-1039, (1974). A nonlinear growth, however, could be considered as a clear signature of the break down of the (finite variance) central limit theorem.
\bibitem{dyn} P. Allegrini et al, Phys. Rev. E., {\bf 54}, 4760-4767, (1996).
\bibitem{barlow} M. T. Barlow, {\em Diffusion on Fractals}, (Springer, 1998).
\bibitem{fgz} U. Frieberg and M. Zahle, Potential Anal. {\bf 16}, 265-277 (2002), U. Freiberg, Forum. Math. {\bf 17}, 87-104, (2005).
\bibitem{pb} Pearson J and Bellissard J,  ArXive: 0802.133v2 (math.OA) (2008).
\bibitem{bukh} Bakhtin Y,  ArXive:0810.0326v2 (math.PR) (2008).
\bibitem{karo} S. Alveberio and W. Karowowski, Stochastic Process Appl. {\bf 53}, 1-22, (1994); J. Math. Phys. {\bf 49}, 093503, (2008).
\bibitem{comment1} Fractional calculus was developed independent of any fractal space. However, more and more studies make it clear that the natural field of applications of FDE are the fractal spaces and stochastic processes(see for intance, \cite{sub, sup}. A search in arXive also reveals the current surge of activity).
\bibitem{comment2} A macroscopic law derived from Lebesgue integration need not be smooth in the sense considered here. A ``fractionally smooth" macroscopic law following a FDE is supposedly nonsmooth (in the ordinary sense) and reflects the presence of randomness.   
\bibitem{dp1}  S. Raut  and D. P. Datta, { Fractals}, {\bf 17}, 45-52 (2009);  ibid, Erratum. {\bf 17}, 547  (2009); ibid, {\bf 18}, (March) (2010), arXive: 1001.1487(math.GM).
\bibitem{dp2} D. P. Datta  and A. Raychaudhuri, { Fractals}, {\bf 18}, June, (2010), arXive:1001.1490 (math.GM), D P Datta and S. Raut, Chaos, Solitons \& Fractals, {\bf 28}, 581-589, (2006).
\bibitem{robin} A. Robinson, {\em Nonstandard Analysis}, North Holland, (Amsterdam, 1966).
\bibitem{dp3} D. P. Datta et al,(2009), Submitted, arXiv:1002.3951v3(math.CA).
\bibitem{dp4} D. P. Datta, { Chaos, Solitons \& Fractals}, {\bf 17}, 621-630, (2003).
\bibitem{comment3} There seems to have been a deep connection between infinitesimals considered here with the Van Hove limit providing a linkage between Hamiltonian dynamics with statistical physics \cite{dyn}, which will be explored elsewhere.
\bibitem{solo} B. Solomyak, Indagationes Math. N. S. {\bf 8}, 133-141, (1997).
\bibitem{note1} If $x\in C_s$ and $t\in C_t$, then $z$, in general, would belong to a Cantor set, when thickness of the original sets satisfy certain restrictions (see, for instance, S. Astels, { Transactions (AMS)}, {\bf 352}, 133-170, (1999)).  
\bibitem{volv} Possibility of numerical fluctuations appears to have been proposed first by I. Volovich, p-adic numbers, ultrametric theory and applications, {\bf 2}, 77-87, (2010).
\end{thebibliography}
\end{document}